\documentclass[10pt]{article}
\usepackage{latexsym}
\usepackage{amsbsy}
\usepackage{amssymb}
\usepackage{amsmath}
\usepackage[latin1]{inputenc}
\usepackage{a4wide}
 \newcommand{\eps}{\varepsilon}
 \newcommand{\R}{\mathbb R}
 \newcommand{\N}{\mathbb N}

 \newcommand{\gs}{\ensuremath{{\mathcal G}} }

 \newcommand{\esm}{\ensuremath{{\mathcal E}_M} }
 
 \newcommand{\ns}{\ensuremath{{\mathcal N}} }

 \newcommand{\comp}{\subset\subset}

 \newcommand{\vphi}{\varphi}

 \newcommand{\D}{{\cal D}}

 \newtheorem{thr}{\hspace*{-3mm} \bf}[section]
 \newcommand{\ethr}{\end{thr}}

 \newcommand{\bd}{\begin{thr} {\bf Definition. }}
 \newcommand{\ed}{\end{thr}}

 \newcommand{\pr}{\noindent{\bf Proof. }}
 
 \newcommand{\ep}{\hspace*{\fill}$\Box$\ms}

 \newcommand{\bthm}{\begin{thr} {\bf Theorem. }}
 \newcommand{\bt}{\begin{thr} {\bf Theorem. }}
 \newcommand{\ethm}{\end{thr}}

 \newcommand{\bp}{\begin{thr} {\bf Proposition. }}
 \newcommand{\bc}{\begin{thr} {\bf Corollary. }}
 \newcommand{\blem}{\begin{thr} {\bf Lemma. }}
 \newcommand{\bex}{\begin{thr} {\bf Example. }\rm} 
 \newcommand{\bexs}{\begin{thr} {\bf Examples. }\rm}
 \newcommand{\brm}{\begin{thr} {\bf Remark. }\rm}

 \newcommand{\beast}{\begin{eqnarray*}}
 \newcommand{\eeast}{\end{eqnarray*}}
 
 \newcommand{\brem}{\begin{thr} {\bf Remark. }\rm}
 \newcommand{\ethi}{\end{thr}}

 \newcommand{\ms}{\medskip\\}
 
 \newcommand{\beq}{ \begin{equation} }
 \newcommand{\eeq}{\end{equation} }
 \newcommand{\bea}{\begin{eqnarray}}
 \newcommand{\eea}{\end{eqnarray}}
 \newcommand{\beas}{\begin{eqnarray*}}
 \newcommand{\eeas}{\end{eqnarray*}}
 \newcommand{\beqs}{\begin{equation*}}
 \newcommand{\eeqs}{\end{equation*}}

 \newcommand{\ben}{\begin{enumerate}}
 \newcommand{\een}{\end{enumerate}}
 \newcommand{\ba}{\begin{array}}
 \newcommand{\ea}{\end{array}}

 \newcommand{\cC}{{\cal C}}
 
 \newcommand{\cG}{{\cal G}}
 \newcommand{\cU}{{\cal U}}
 
 \newcommand{\cN}{{\cal N}}
 \newcommand{\cE}{{\cal E}}
 
 \newcommand{\cP}{{\cal P}}

 \newcommand{\cbesk}{{\cal C}^{\infty}}
 
 \newcommand{\rp}{{\mathbb R}^p}
 \newcommand{\rqq}{{\mathbb R}^q}
 \newcommand{\rn}{{\mathbb R}^n}
 \newcommand{\rmm}{{\mathbb R}^m}

 \newcommand{\rl}{{\mathbb R}^l}

 \newcommand{\vf}{{\bf v}}
 \newcommand{\pd}{\partial}
 
 \newcommand{\un}{u^{(n)}}
 \newcommand{\rtilda}{\widetilde{\R}}
 \newcommand{\prn}{{\rm pr}^{(n)}}

 \newcommand{\deta}{\frac{d}{d \eta}{\Big\vert}_{_{0}}}

\begin{document}

 

\title{ Generalized Group Actions in a Global Setting}
\author{Sanja Konjik \footnote{Faculty of Agriculture, University of Novi Sad, Trg Dositeja Obradovi\' ca 8, 21000 Novi Sad, SCG.
Electronic mail: kinjoki@neobee.net}\\
         Michael Kunzinger \footnote{Faculty of Mathematics, University of Vienna, Nordbergstr.\ 15, A-1090 Wien, Austria,
         Electronic mail: michael.kunzinger@univie.ac.at}\\ 
       }

\date{}
\maketitle
                                                                           
\begin{abstract}
We study generalized group actions on differentiable manifolds in the Colombeau framework,
extending previous work on flows of generalized vector fields and symmetry group analysis
of generalized solutions. As an application, we analyze group invariant generalized functions
in this setting.

\end{abstract}

\section{Introduction} \label{intr}
Lie group analysis of differential equations is an indispensable tool for
studying invariance properties of solutions of PDE as well as for finding
explicit solutions, with a wealth of applications (cf.\ \cite{BK, Olv}).
In \cite{Me, S, SZ, Zie}, a study of invariance properties of distributions
and distributional
solutions of linear partial differential equations was initiated. 
Later on, symmetry group analysis of PDEs in generalized functions and 
systematic  methods of deriving group invariant fundamental 
solutions using infinitesimal techniques of group analysis were developed
\cite{ga, wi, bi, amp}. Clearly, in the distributional setting a restriction
to linear equations and linear projectable transformation groups 
is unavoidable. On the other hand, many applied problems (e.g.,
systems of conservation laws) underline the need for an extension of the
above techniques in order to handle  nonlinear problems
involving singularities. Algebras of generalized functions provide a setting
for addressing such questions in a coherent way. This line of research was 
initiated in \cite{OR, RW, RW1} in the framework of the 'nowhere dense' algebras
of E.E.\ Rosinger. 

An alternative approach, based on Colombeau's theory of algebras of generalized
functions (\cite{c1,c2,MOBook}), was developed in \cite{symm, DKP, ObCont, MoRot}
and will form the basis for the present paper.
In particular, in \cite{symm, DKP} criteria for classical symmetry 
groups to transform weak (distributional, Colombeau or associated)
solutions of a given (smooth) system of differential equations
into other solutions of the same type were given. In \cite{symm, DKP, 
ObCont, MoRot},  additionally both the differential operators and the group actions 
are allowed to be given by generalized functions. The setting of generalized functions employed
in these works is that of $\cG_{\tau}$, the space of tempered Colombeau functions.
As elements of $\cG_{\tau}$ are characterized by global bounds, this setting appears
unsuitable for an extension of the theory to the manifold setting. To lift this
limitation, in the present work we employ the recently developed theory of
Colombeau generalized functions taking values in differentiable manifolds
(\cite{gfvm, gfvm2}) as well as the theory of generalized flows of singular
vector fields (\cite{flows}) to extend symmetry group analysis in Colombeau
generalized functions to a global setting.

The paper is divided into 6 sections. Section 1 
provides  basic notations and definitions from Colombeau's
theory of algebras of generalized functions (in particular in the manifold
setting) and symmetry group analysis.
In section \ref{notations} we consider generalized group actions and provide a notion
of  rank of a generalized function, which will be crucial for the infinitesimal
criteria to be developed in section \ref{infcriter}.  The question of localizing 
Colombeau generalized functions and an analysis of solution sets of generalized
equations is the focus of section \ref{localization}. By borrowing a notion from nonstandard analysis
we introduce the concept of near-standard points and show that these
suffice to characterize equality of Colombeau functions.
In section \ref{infcriter} we prove an infinitesimal criterion for symmetry groups of
generalized algebraic equations and apply the obtained results
in section 5 to symmetry group analysis of differential equations
in the Colombeau framework. Finally, in section \ref{rotsect} we turn to the
topic of group invariant generalized functions in this setting. 
Based on a recent result of Pilipovi\'{c}, Scarpalezos and Valmorin (\cite{PSV})
we provide an affirmative  answer to an open question posed by M.\ Oberguggenberger
in  \cite{MoRot} whether standard  rotations suffice to characterize rotational invariance of
Colombeau generalized functions. 

\section{Notations}\label{notations}
 In what follows, $M$
 and $N$ will denote smooth, connected, paracompact Hausdorff
 manifolds of dimensions $m$ and $n$, respectively.

 Set $I=(0,1]$ and denote by ${\cal P}(M)$ the space of linear
 differential operators on $M$. The spaces of moderate resp.\ negligible nets in $M$ are
 defined as 
 \beas
 {\cal E}_M(M)&:=&\{ (u_\eps)_\eps\in{\mathcal C}^\infty(M)^I:\ \forall
 K\subset\subset M,\ \forall P\in{\cal P}(M)\ \exists p\in\N:\\
 &&\hphantom{(u_\eps)_\eps\in{\mathcal E}(X):\ \forall
  K\subset\subset X,\ \forall P\in{\cal P}(\}} \sup_{x\in
  K}|Pu_\eps(x)|=O(\eps^{-p})\}\\
  {\cal N}(M)&:=& \{ (u_\eps)_\eps\in\cE_M(M):\ \forall K\subset\subset M,\ 
 \forall m \in\N_0:\ \sup_{x\in K}|u_\eps(x)|=O(\eps^{m})\}\, 
 \eeas
 (due to \cite{book}, ch.\ 1, Th.\ 1.2.3, for the characterization of
 $\cN(M)$ as a subspace of $\esm(M)$ it is sufficient to estimate only 
 the $0$-th order derivative).
 Clearly, $\cN(M)$ is an ideal of the differential algebra $\cE_M(M)$. 
 The special Colombeau algebra $\cG(M)$ on $M$ is defined as the quotient
 $\cE_M(M)/\cN(M)$; it is an associative, commutative differential algebra whose
 elements are equivalence classes denoted by $u=[(u_\eps)_\eps]$.
 $\cG(\_)$ is a fine sheaf of differential algebras with respect to the
 Lie derivative along smooth vector fields. $\cC^\infty(M)$ is a subalgebra
 of $\cG(M)$ and there exist injective sheaf morphisms embedding $\D'(\_)$ linearly 
 into $\cG(\_)$.

 A point value characterization of Colombeau generalized functions
 is based on the concept of compactly supported generalized
 points (\cite{point, found}). The space of compactly supported generalized points $M_c$ is
 the set of all nets $(x_\eps)_\eps\in M^I$ for which $x_\eps$
 stays in a fixed compact set for $\eps$ small. In $M_c$ one introduces
 an equivalence relation $\sim$ in the following way: for 
 $(x_\eps)_\eps,\, (y_\eps)_\eps\in M_c$, 
 $(x_\eps)_\eps\sim (y_\eps)_\eps \Leftrightarrow 
 d_h(x_\eps,y_\eps)=O(\eps^m)$,
 for each $m>0$, where $d_h$ denotes the distance function induced on
 $M$ by one (hence any) Riemannian metric $h$. The quotient space
 $\widetilde{M}_c:=M_c/\sim$ is called the space of compactly
 supported generalized points on $M$, and we denote its elements 
 by $\tilde{x}=[(x_\eps)_\eps]$. In the case $M=\R$
 one also defines the ring of generalized numbers
 $\widetilde{\R}$ as the quotient of the set
 of moderate nets of numbers $(r_\eps)_\eps \in \R^I$ with 
 $|r_\eps| = O(\eps^{-p})$ for some $p\in\N$ modulo the set of
 negligible nets $(r_\eps)_\eps$ with $|r_\eps| = O(\eps^{m})$ for each $m$. It is the ring of
 constants in the Colombeau algebra. Insertion of a compactly supported
 generalized point into any representative of a Colombeau generalized
 function produces a well-defined element of $\widetilde{\R}$. Moreover,
 elements of $\cG(M)$ are uniquely determined by their values on
 $\widetilde{M}_c$.

 In order to describe generalized functions on the manifold
 $M$ taking values in the manifold $N$ one introduces the
 space $\cG[M,N]$ of compactly supported (or c-bounded for short) 
 generalized functions. A net $(u_\eps)_\eps\in \cC^\infty(M,N)^I$
 is called c-bounded if
 $$
 \forall K\comp M\ \exists \eps_0>0\ \exists K'\comp N\ \forall
 \eps<\eps_0:\ u_\eps(K) \subseteq K'.
 $$ 
 A c-bounded net is moderate if it satisfies:
 \begin{itemize}
 \item[] $\forall k\in\N$, for each chart $(V,\vphi)$ in $M$, each 
   chart $(W,\psi)$ in $N$, each $L\comp V$ and each $L'\comp W$
   there exists $p\in \N$  with
   $$\sup\limits_{x\in L\cap u_\eps^{-1}(L')} \|D^{(k)}
   (\psi\circ u_\eps \circ \vphi^{-1})(\vphi(x))\| =O(\eps^{-p})\,,
   $$
 \end{itemize}
 Denote by $\cE_M[M,N]$ the set of all moderate c-bounded nets. Introduce
 an equivalence relation $\sim$ in $\cE_M[M,N]$ in the following way: 
 $(u_\eps)_\eps, (v_\eps)_\eps \in \cE_M[M,N]$, 
 $(u_\eps)_\eps \sim (v_\eps)_\eps$ if
 \begin{itemize}
 \item[(i)] $\forall K\comp M$, $\sup_{x\in K}d_h(u_\eps(x),v_\eps(x)) \to 0$
 ($\eps\to 0$)
 for some (hence every) Riemannian metric $h$ on $N$.
 \item[(ii)] $\forall k\in \N_0\ \forall m\in \N$,
 for each chart    $(V,\vphi)$
   in $M$, each chart $(W,\psi)$ in $N$, each $L\comp V$
   and each $L'\comp W$:
 $$
 \sup\limits_{x\in L\cap u_\eps^{-1}(L')\cap v_\eps^{-1}(L')}\!\!\!\!\!\!\!\!\!\!\!\!\!\!\!\!\!\!\!
 \|D^{(k)}(\psi\circ u_\eps\circ \vphi^{-1}
 - \psi\circ v_\eps\circ \vphi^{-1})(\vphi(x))\|
 =O(\eps^m).
 $$
 \end{itemize}
 The space of c-bounded Colombeau generalized functions from $M$ to $N$
 is defined as the quotient $\cG[M,N]:=\cE_M[M,N]/\sim$.

 Alternative characterizations
 of the notions of moderatness and equivalence for the elements of
 $\cC^\infty(M,N)^I$ are: $(u_\eps)_\eps\in\cE_M[M,N]\Leftrightarrow
 (f\circ u_\eps)_\eps\in \cE_M(M)$, $\forall f\in\cC^\infty(N)$
 (\cite{gfvm2},\ Prop.\ 3.2)
 and $(u_\eps)_\eps \sim (v_\eps)_\eps\ ((u_\eps)_\eps, (v_\eps)_\eps)
 \in \cE_M[M,N]) \Leftrightarrow (f\circ u_\eps -
 f\circ v_\eps)_\eps\in \cN(M)$, $\forall f\in\cC^\infty(N)$
 (\cite{gfvm2},\ Th.\ 3.3).
 Similarly as for the elements of $\cG(M)$, if $u\in\cG[M,N]$ and
 $\tilde{x}\in\widetilde{M}_c$ then $u(\tilde{x})=[(u_\eps(x_\eps))_\eps]$
 is a well-defined element of $\widetilde{N}_c$, and
 elements of $\cG[M,N]$ are uniquely determined by their values
 on all compactly supported generalized points on $M$, i.e.
 $u=v \Leftrightarrow u(\tilde{x})=v(\tilde{x})$, 
 $\forall \tilde{x}\in\widetilde{M}_c$ (\cite{gfvm2},\ Th.\ 3.5).

 If $E \to M$ is any vector bundle over $M$, denote by
 $\Gamma(M,E)$ the space of smooth sections of $E$, and
 by $\cP(M,E)$ the space of differential operators 
 $\Gamma(M,E)\to\Gamma(M,E)$. The module of
 generalized sections of $E$, $\Gamma_{\cG}(M,E)$, is defined
 as the quotient $\Gamma_{\cE_M}(M,E)/\Gamma_{\cN}(M,E)$ where
 \beas
 \Gamma_{\cE_M}(M,E)&:=& \{ (s_\eps)_{\eps}\in \Gamma(M,E)^I :
       \ \forall P\in {\mathcal P}(M,E)\, \forall K\comp M \, \exists p\in \N:\\
 &&\hspace{6cm}\sup_{x\in K}\|Pu_\eps(x)\|_h = O(\eps^{-p})\}\\
 \Gamma_{\cN}(M,E)&:=& \{ (s_\eps)_{\eps}\in \Gamma_{\cE_M}(M,E) :
       \ \forall K\comp M \, \forall m\in \N:\\
 &&\hspace{6.5cm}\sup_{x\in K}\|u_\eps(x)\|_h = O(\eps^{m})\},
 \eeas
 where $\|\,\|_h$ is the norm on the fibers of $E$ induced by
 any Riemannian metric on $M$. 
 $\Gamma_{\cG}(\_\,,E)$ is a fine sheaf of  projective and finitely generated  
 $\cG(M)$-modules, and
 $$ 
 \Gamma_{\cG}(M,E)=\cG(M)\otimes_{\cC^\infty(M)}\Gamma(M,E).
 $$
 If $E$ is some tensor bundle $T^r_sM$ we write $\cG^r_s(M)$
 instead of $\Gamma_\cG(M,T^r_sM)$; in particular, if $E$ is
 the tangent bundle $TM(=T^1_0M)$ then $\cG^1_0(M)$ is the
 space of generalized vector fields on $M$.

 We say that a generalized vector field $\xi\in \cG^1_0(M)$ is
 locally bounded resp.\ locally of $L^\infty$-log-type 
 if for all $K\comp M$ and one (hence every) Riemannian metric
 $h$ on $M$ we have for any representative $(\xi_\eps)_\eps$
 and $\eps$ sufficiently small 
 $$
 \sup_{x\in K}\|\,\xi_\eps|_x\,\|_h\,\leq\,C\, \quad \mbox{resp.} 
 \quad \sup_{x\in K}\|\,\xi_\eps|_x\,\|_h\,\le C |\log \eps|\,.
 $$
  $\xi$ is called globally bounded with respect to $h$ if for
 some (hence every) representative $(\xi_\eps)_\eps$ of $\xi$
 there exists $C>0$ with
 $$
 \sup_{x\in M}\|\,\xi_\eps|_x\,\|_h\,\leq\, C ,
 $$
 for $\eps$ small (cf.\ \cite{flows}, Def.\ 3.4).

 To conclude this section we fix some notations from symmetry group analysis
 of differential equations, following \cite{Olv}.
 Let $X$ and $U$ be spaces of  independent and dependent variables and 
 suppose that $G$ is a local
 Lie group of transformations acting regularly on some open subset
 $M\subseteq X\times U$; for the group action we write $g\cdot(x,u)=
 (\Xi_g(x,u), \Psi_g(x,u))$, with appropriate smooth functions
 $\Xi_g$ and $\Psi_g$. If $\Xi_g$ does not depend on the
 dependent variables the group action is called projectable.
 The $n$-jet space of $M$ will be
 denoted by $M^{(n)}$ and the $n$-th prolongation of a group
 action $g$ resp.\ vector field $\vf$ by $\prn g$ resp. 
 $\prn \vf$.
 If $\Delta_\nu(x, \un)=0$ $(1\leq \nu\leq l)$ is a system
 of $n$-th order differential equations on $M$, where $\Delta:
 X\times U^{(n)}\to \R^l$ is a smooth function, then the
 solution set of $\Delta$ is the subvariety $S_\Delta:=\{
 (x, \un):\ \Delta(x,\un)=0\}$ of $X\times U^{(n)}$. We say
 that a function $f$ is a solution of the system if the
 $n$-jet of the graph $\Gamma_f=\{(x,f(x)):\ 
 x\in\Omega\}\subset X\times U\}$ of $f$, i.e. 
 $\Gamma^{(n)}_f$ is contained in $S_\Delta$. A symmetry group
 of $\Delta$ is a local transformation group $G$ acting on
 $M$ with the property that whenever $u=f(x)$ is a solution
 of the system and $g\cdot f$ $(g\in G)$ is defined, then
 $g\cdot f$ is again a solution of $\Delta$.

\section{Generalized Group Actions} \label{gga}

 To begin with we recall the following definitions from \cite{flows}:

 \bd \label{gen action}
 A generalized group action on a manifold $M$ is
 an element $\Phi \in \cG[\R\times M,M]$ with the
 following properties:
 \begin{itemize}
 \item[(i)] $\Phi(0, \cdot) = \mbox{id} \;\;$ in $\cG[M,M]$
 \item[(ii)] $\Phi(\eta_1+\eta_2, x)=\Phi(\eta_1,\Phi(\eta_2, x)) 
  \;\;$ in $\cG[\R^2\times M, M]$.
 \end{itemize}
 \ethi

In the following definition we make use of 
$\gs^h$, the space of hybrid Colombeau functions
defined on a manifold and taking values in a vector bundle which was
introduced in \cite{gprg} (see also \cite{book}).
 \bd
 Let $\xi \in \cG^1_0(M)$ be a generalized vector field  such that
 there exists a unique generalized group action  $\Phi\in\cG[\R\times 
M,M]$ satisfying
 \beq \label{infgen}
 \frac{d}{d\eta}\Phi(\eta,x) = \xi(\Phi(\eta,x)) \qquad \mbox{in }
 \cG^h[\R\times M, TM]
 \eeq
 Then $\xi$ is called the infinitesimal generator of $\Phi$ and both 
 $\xi$ and its generalized flow $\Phi$ are called $\cG$-complete.
 We call $\xi$ and $\Phi$ strictly $\gs$-complete if, in addition,
there exist representatives $(\xi_\eps)_\eps$, $(\Phi_\eps)_\eps$
such that $\Phi_\eps$ is the flow of $\xi_\eps$ for each $\eps \in I$.
 \ethi

Even for not necessarily $\gs$-complete group actions $\Phi$ we shall 
call a generalized vector field $\xi$ an infinitesimal generator of
$\Phi$ if (\ref{infgen}) holds. In practice, since in order to show
$\gs$-completeness one usually works componentwise, the condition of
strict $\gs$-completeness is normally no additional restriction, cf.\
the following remark.

 \brem
 Sufficient conditions for $\gs$-completeness of a generalized 
 vector field $\xi$ have been derived in \cite{flows}, Th.\ 3.5
 for the case of $(M,h)$ a complete Riemannian manifold, to wit:
 \begin{itemize}
 \item[(i)] $\xi$ is globally bounded with respect to $h$, and 
 \item[(ii)] for each first-order differential operator $P\in
 {\cal P}(M,TM)$, $P\xi$ is locally of $L^\infty$-log-type.
 \end{itemize}
In fact, these conditions even ensure strict $\gs$-completeness of 
$\xi$.
 \ethi

One of our main interest in generalized group actions in this work will
be symmetry properties in the following sense: 
 \bd\label{algsymdef}
 Let $F\in\cG(M)$ and let $\Phi$ be a $\gs$-complete generalized group
 action on $M$. $\Phi$ is called a symmetry group of the
 equation
 $$
 F(x)=0
 $$
 in $\cG(M)$ if for any $\tilde{x} \in \widetilde{M}_c$
 with $F(\tilde{x})=0 \in \widetilde \R$ we have $F(\Phi(
 \tilde{\eta}, \tilde{x}))=0$ in $\widetilde \R$, for every
 $\tilde{\eta} \in \widetilde{\R}_c$ (i.e., $\eta
 \mapsto F(\Phi(\eta,\tilde{x}))=0$ 
 in $\cG(\R)$). If $F= (F_\nu)_{\nu=1}^l\in \gs(M)^l$, then $\Phi$ is 
 called a symmetry group of the equation $F=0$ if it is a symmetry
 group of each equation $F_\nu=0$ $(1\le \nu\le l)$.
 \ethi

 We note that, since $\gs[M,\R]$ is naturally contained in $\gs(M)$,
 the above definitions and results directly apply to c-bounded 
 generalized functions as well.

 As in the classical case (cf.\ \cite{Olv}, ch.\ 2) our first aim is to
 derive infinitesimal criteria characterizing symmetries of 
 ``algebraic'' equations as in \ref{algsymdef}. In the smooth setting,
 one supposes a maximal rank condition on $F$ and then uses
 distinguished local charts to obtain the desired result. In our 
present context,
 however, a direct transfer of classical
 methods is impossible due to the lack of structure of the
 space $\widetilde{M}_c$ of compactly supported generalized points
 on $M$. In particular, elements of $\widetilde{M}_c$ are only very 
 weakly localized in the sense that every $\tilde{x} \in 
\widetilde{M}_c$ 
 possesses a representative contained in a suitable compact set
 in $M$. We therefore call an open set $U\subseteq M$
 a neighborhood of $\tilde{x} = [(x_\eps)_\eps]$ if 
 $$
 \exists \eps_0\ \exists K\subset\subset U\ \mbox{s.t. }
 x_\eps\in K \ \forall \eps<\eps_0\,.
 $$ 
 Moreover, in the absence of an inverse function theorem in the
 generalized function setting, it is a priori not clear how to define
 the rank of a generalized function. Since, on the positive side, 
inversion 
of generalized functions is possible in $\gs[M,N]$
 we suggest the following notion of rank of a generalized map:

 \bd \label{def rank}
Let $F\in \gs[M,N]$, $k\in 
\{0,\dots,\min (m,n)\}$ and $\tilde{x} \in \widetilde{M}_c$.
$F$ is called of rank $k$ in $\tilde{x}$ if there exist
open neighborhoods $U\subseteq M$ of $\tilde{x}$, $V\subseteq N$
of $F(\tilde{x})$, open sets $U'\subseteq \R^m$, 
$V'\subseteq \R^n$, $\eps_0>0$ and diffeomorphisms $\vphi_\eps: U\to 
U'$,
$\psi_\eps: V\to V'$ for each $\eps\in (0,\eps_0]$ with  $\vphi = 
[(\vphi_\eps)_\eps]
\in \gs[U,U']$, $\vphi^{-1}:=[(\vphi_\eps^{-1})_\eps] \in \gs[U',U]$, 
$\psi = [(\psi_\eps)_\eps]\in \gs[V,V']$
$\psi^{-1}:=[(\psi_\eps^{-1})_\eps] \in \gs[V',V]$ such that 
$F|_U\in \gs[U,V]$ and
$$
\psi\circ F\circ \vphi^{-1} = (x_1,\dots,x_m) \mapsto 
(x_1,\dots,x_k,0,\dots,0)
$$
in $\gs[U',V']$. If $A\subseteq U$ then $F$ is called of rank $k$ 
globally
on $A$.
\ethi

It is straightforward to adapt this definition also to the case
where $F\in \cG(M)^l$ (set $N=\rn$ and $\psi_\eps=\mathrm{id}$ for all 
$\eps$).

According to the above discussion  it is natural to ask whether a more 
strict localization 
than the one used in Definition \ref{def rank} is attainable in 
general. 
Before we proceed with the theory of symmetry groups of generalized 
algebraic equations we should therefore investigate the possibility of 
localizing
Colombeau generalized functions resp.\ solution sets of generalized
equations. The following section is devoted to this purpose.

\section{Localization} \label{localization}

By the point value characterization of Colombeau generalized functions 
(cf.\ \cite{point}, \cite{ndg}, \cite{gfvm2}), elements of $\cG(M)$ as 
well as of $\cG[M, N]$ are
uniquely determined by their values on compactly supported
generalized points on $M$. 

As was mentioned in the previous section, elements of
$\widetilde{M}_c$ are only weakly localized, so in particular the 
existence
of suitable open neighborhoods of $\tilde{x}\in \widetilde{M}_c$ as
in Definition \ref{def rank} is not necessarily guaranteed.
Therefore the question arises whether we need all elements of
$\widetilde{M}_c$ to characterize elements of $\cG(M)$ (or $\cG[M, 
\R]$)
or if more strongly localized generalized points suffice. The following
definition borrows a concept from nonstandard analysis to specify what
is meant by this notion:

\bd
A point $\tilde{x}\in \widetilde{M}_c$ is called near-standard
if there exists $x\in M$ such that $\tilde{x} \approx x$ (i.e.,
$x_\eps \to x$ $(\eps\to 0)$ for every representative of $\tilde{x}$).
\ethi
In particular,  any neighborhood of $x$ is a neighborhood of $\tilde{x} 
\approx x$ 
in the sense of section \ref{gga}. Near-standard points indeed suffice 
to 
characterize Colombeau generalized functions:

 \bp
 (i) Let $u\in\cG(M)$. Then $u=0$ if and only if
 $u(\tilde{x})=0$, for all near-standard points
 $\tilde{x} \in \widetilde{M}_c$.
 
 (ii) Let $u, v \in \cG[M,N]$. Then $u=v$ if and only if
 $u(\tilde{x})=v(\tilde{x})$, for all near-standard points
 $\tilde{x} \in \widetilde{M}_c$.  
 \ethi

 \pr
 (i) One direction is clear. So, suppose that 
 $u(\tilde{x})=0$, for all near-standard points
 $\tilde{x} \in \widetilde{M}_c$
 and suppose that $u\not= 0$. Then
 \beq \label{near-standard}
 \exists K\subset\subset M\ \exists m\ \forall k\in \N\
 \exists x_k\in K\ \exists \eps_k<\min(\frac{1}{k},\eps_{k-1})\ \quad
 |u_{\eps_k}(x_k)|>k\eps_k^m.
 \eeq
 Since $K$ is a compact set there exists a subsequence
 $x_{k_l}$ which converges to $x\in K$. Set
 $x_\eps:=x_{k_l}$ for $\eps \in (\eps_{k_{l+1}}, \eps_{k_l}]$
 and $\tilde{x} := [(x_\eps)_\eps]\in \widetilde{M}_c$.
 $\tilde{x}$ is a near-standard point and from
 (\ref{near-standard}) it follows that $u(\tilde{x}) \not= 0$,
 which gives a contradiction.

 (ii) Necessity is again obvious. For the converse direction we use the
 characterization of $c$-bounded generalized functions given in
 \cite{gfvm2}. Let $f\in \cbesk(N)$. Then $f\circ u$ and $f\circ v$
 are well defined elements of $\cG(M)$. For any near-standard point
 $\tilde{x}\in \widetilde{M}_c$ we have by (i) that
 $$
 (f \circ u)(\tilde{x}) = (f\circ v)(\tilde{x})\,, 
 $$
so  $f\circ u=f\circ v$ in $\gs(M)$.
 Hence, $(f\circ u_\eps - f\circ v_\eps)_\eps \in \cN(M)$ and by 
 \cite{gfvm2}, Th. 3.3 it follows that $u=v$.   
 \ep
 
In the smooth setting, a maximal rank
condition on the set of solutions of an equation $F(x)=0$ allows to
derive an infinitesimal criterion for symmetry groups of the equation
(cf.\ \cite{Olv}, ch.\ 2). In the generalized case, however,
the assumption of maximal rank in each near-standard point
$\tilde{x}\in \widetilde{M}_c$ which is a solution of
$F(x)=0$, $F\in \cG(M)$, may be insufficient. We illustrate
this by the following example:

 \bex \label{counter ex}
 Set $M=\R$, $I=(0,1]$ and
 $J=\bigcup_{n=1}^\infty (\frac{1}{2n+1}, \frac{1}{2n}]$.
 Let
 \begin{equation} \label{solution}
 x_\eps = \left\{\begin{array}{rl}
                  0, & \eps \in J,\\
                  1, & \eps \in I\setminus J
                  \end{array}
           \right.
 \end{equation}
 \begin{equation} \label{F}
 F_\eps(x) = \left\{\begin{array}{rl}
                  x, & \eps \in J,\\
                  x -1, & \eps \in I\setminus J.
                         \end{array}
                  \right.
 \end{equation}
 Then $F_\eps(x_\eps)=0$ for all $\eps$ and 
$\tilde{x}:=[(x_\eps)_\eps]$ is
 not a near-standard point. We claim that the solution set 
$S_F=\{\tilde{y} \in \widetilde{\R}_c\mid 
F(\tilde{y})=0\}$ does not contain any
 near-standard point. To see this, suppose that
 $\tilde{y}=[(y_\eps)_\eps]$ satisfies $F(\tilde{y})=0$.
 Then
 $$
 F_\eps(y_\eps)= y_\eps -1 \qquad \mbox{ on }\ I\setminus J
 $$ 
 and
 $$
 F_\eps(y_\eps)=y_\eps \qquad \mbox{ on }\ J.
 $$
Suppose that $\tilde{y}$ is a near-standard point and choose
$y\in \R$ such that $y_\eps \to y$ when $\eps \to 0$.
Then since $F(\tilde y) = 0$ we obtain $y-1=0=y$, a contradiction. Moreover, the above reasoning 
implies that $\tilde{x}$ is in fact the only zero of the
 equation $F(x)=0$ in $\widetilde{\R}_c$.
 \ethm

 This example shows that there exist functions whose solution
 set is nonempty although it does not contain any near-standard points.
 In order to obtain infinitesimal criteria for an
 equation $F(x)=0$ we will therefore have to require a maximal rank 
condition in
a neighborhood of all of $S_F$, no matter which types of generalized 
 points belong to it. 

An alternative localization strategy consists in considering an open
covering $\cU$ of $M$.
Since $\gs$ is a sheaf, $F=0$ on $M$ if and only if $F=0$ on each open
 set $U\subseteq M$. Also, if $\Phi$ is a symmetry
 group of the equation
 $$
 F(x)=0 \mbox{ in } \cG(M)
 $$
 then for every open covering $\cU$ of $M$, $\Phi$ is a symmetry
 group of
 $$
 F|_U (x)=0
 $$
for each $U\in \cU$.
However, a localization to near-standard points fails in general:
consider again Example \ref{counter ex}. Let 
$U_1 =(-\infty, \frac{1}{2})$ and $U_2 =(\frac{1}{4}, \infty)$. Then
the intersection of $S_F$ with both
$(\widetilde{U}_1)_c$ and $(\widetilde{U}_2)_c$ is empty, although
$S_F$ itself is nonempty, consisting precisely of the generalized 
point $\tilde{x}$ from (\ref{solution}).

\section{Infinitesimal criteria} \label{infcriter}
Our aim in this section is to derive infinitesimal criteria for 
symmetry groups
of algebraic equations in the Colombeau setting. To this end we will 
need the following auxilliary result:

\blem \label{help1}
Let $\psi_\eps: M \to N$ $(\eps\in I)$ be a net of diffeomorphisms 
such that
$\psi=[(\psi_\eps)_\eps] \in \gs[M,N]$ and 
$\psi^{-1}=[(\psi_\eps^{-1})_\eps] 
\in \gs[N,M]$. If $\Phi$ is a strictly $\gs$-complete group action on 
$N$ with
generator $\xi\in \gs^1_0(N)$ then $\psi^* \Phi = 
[(\psi^*\Phi_\eps)_\eps]$
is a stricly $\gs$-complete group action on $M$ with infinitesimal 
generator
$\psi^* \xi = [(\psi^* \xi_\eps)_\eps]$. 
\ethi
\pr Choose representatives $(\xi_\eps)_\eps$, $(\Phi_\eps)_\eps$ as in 
the
definition of strict $\gs$-completeness. Then for each fixed $\eps\in 
I$, 
$\psi^*\Phi_\eps(\eta,x) = \psi^{-1}\circ 
\Phi_\eps(\eta,\psi(x))$ is a group action on $M$ with generator 
$\psi^*\xi_\eps
= T\psi^{-1}\circ \xi_\eps \circ \psi$. Since equation (\ref{infgen}) 
transfers
componentwise from $N$ to $M$, strict $\gs$-completeness of the  
pullback follows.
\ep

In the formulation of Theorem \ref{glavna} below 
we will make use of the following definition:
 a subset of $\R^n$ is called an $n$-dimensional box if it is
 a product $I_1\times \dots \times I_n$ of $n$ finite or
 infinite open intervals in $\R$.
\bt \label{glavna}
Let $\Phi$ be a stricly $\cG$-complete group action on $M$ with 
generator $\xi$. 
Let $F\in \cG(M)^l$ be of maximal rank on some $U$ with
 $\widetilde{U}_c\supseteq \Phi((-\eta_0, 
\eta_0)^\sim_c \times S_F)$ ($\eta_0>0$) via a generalized chart 
$\psi\in\gs[U,V]$, 
where $S_F:=\{\tilde{x} \in \widetilde{M}_c | F(\tilde{x})=0\}$.
Set $\bar \xi := (\psi^{-1})^*\xi$ and suppose that one of the 
following 
conditions holds:
\begin{itemize}
\item[(i)] $V=\R^m$ and $\bar \xi$ possesses a representative $(\bar 
\xi_\eps)_\eps$ 
satisfying: 
$$
\exists C,\, \eps_0 > 0 \ \mbox{ such that } 
|\bar\xi_\eps(x)|\le C(1+|x|) \qquad (x\in\R^n,\ \eps<\eps_0).
$$
\item[(ii)]  
$V$ is a box and $D\bar\xi$ is locally of $L^\infty$-$\log$-type. 
\end{itemize}
Then $\Phi$ is a symmetry group of
\beq \label{system F}
F_\nu(x)=0 \qquad\qquad (1\leq \nu \leq l)
\eeq
if and only if 
\beq \label{inf criterion alg}
\xi(F_\nu)|_{\tilde{x}} = 0 \qquad 1\leq \nu \leq l \quad
\forall \tilde{x}\in S_F.
\eeq
\ethi

\pr
Let $\Phi$ be a symmetry group of (\ref{system F}). Then for
each $\tilde{x} \in S_F$, the generalized function 
 $$
 \eta \mapsto F(\Phi(\eta, \tilde{x})) 
 $$
equals $0$ in $\gs(\R)$. Therefore,
 $$
 0=\left.\frac{d}{d\eta}\right\vert_{0}F(\Phi(\eta,\tilde{x}))
 =\xi(F)|_{\tilde{x}}.
 $$ 
Conversely, by assumption we have
$F\circ \psi^{-1}=\mbox{pr}: V\subseteq  \rmm \to \rl$. By Lemma 
\ref{help1}, 
 $(\psi^{-1})^*\xi$ is a strictly $\cG$-complete vector field
 on $V$ with flow $\bar \Phi (\eta,x):=(\eta, x)
 \mapsto \psi \circ \Phi(\eta, \psi^{-1}(x))$.
 Write
 $$\bar \xi=
 (\psi^{-1})^*\xi = \sum_{i=1}^m \bar \xi_i \pd_{x_i}
 $$
 and
 $$
\bar F:= (\psi^{-1})^*F = (x_1, \dots, x_m)
 \mapsto (x_1, \dots, x_l).
 $$
 Then
 $
 S_{\bar F}=\psi(S_F)=\{\tilde{x}\in \widetilde{V}_c |\,
 \bar F(\tilde{x})=0 \mbox{ in } \widetilde{\R}^l_c\}=
 \{\tilde{x} \in \widetilde{V}_c|\, (\tilde{x}_1, \dots,
 \tilde{x}_l)=0 \mbox{ in }\widetilde{\R}^l_c\}.
 $
 Moreover, $\bar \xi(\bar F)=0$ on $S_{\bar F}$
 means that $\bar \xi_i|_{V\cap (\{0\}\times \R^{m-l})}=0$ in
 $\cG(V\cap (\{0\}\times \R^{m-l}))$, for $i=1, 2,\dots , l$.
 Hence, $\bar \xi|_{V\cap (\{0\}\times \R^{m-l})}$ has
 a representative $(\bar \xi_\eps)_\eps$ with $\bar \xi_{1\eps},
 \dots, \bar\xi_{l\eps}\equiv 0$. Write $\bar \xi=
 (\bar \xi', \bar \xi'')\in \cG(V)^l\times
 \cG(V)^{m-l}$ and let $\tilde{x}\in S_{\bar F}$.
 Then $\tilde{x}$ has a representative $(x_\eps)_\eps$
 such that $x_\eps=(0,x_\eps'')\in (\rl\times \R^{m-l})\cap V$
 for all $\eps$. 

Suppose now that assumption (i) is satisfied. Then the initial value 
problem
\begin{equation} \label{IVP1}
\displaystyle
  \begin{array}{rcl}
 \frac{d}{d\eta}\phi(\eta) & = & \bar \xi''(0,\phi(\eta))\\[5pt]
 \phi(0) & = & \tilde{x}'',
  \end{array}
\end{equation}
possesses a solution on $\{0\}\times \R^{m-l}$ (see the existence part 
of
the proof of Th.\ 3.2 in \cite{flows}).  
Set $\tilde{x}'':=[(x_\eps'')_\eps]$ with $x_\eps''$ as above. 
Let $\phi$ be a solution of (\ref{IVP1}).
Then $\bar \Phi(\eta, \tilde{x})= (0,\phi(\eta))$.
Indeed, let $\widetilde{\Phi}(\eta, x):= (0, \phi(\eta))$. Then
$$
\widetilde{\Phi}(0, \tilde{x})= (0, \tilde{x}'')=\tilde{x}
$$
and
$$
\frac{d}{d\eta}\widetilde{\Phi}(\eta, \tilde{x}) =
(0, \phi'(\eta))=(0,  \bar\xi''(0, \phi(\eta)))=
\bar\xi(0, \phi(\eta))= \bar\xi (\widetilde{\Phi}
(\eta, \tilde{x})).
$$
Hence $\bar \Phi$ and $\widetilde{\Phi}$ both solve
the initial value problem 

\begin{equation}\label{IVP2}
\begin{array}{rcl}
 \frac{d}{d\eta}\Phi(\eta) & = & \bar\xi(\Phi(\eta))\\[5pt]
 \Phi(0) & = & \tilde{x}.
\end{array}
\end{equation}

Since $\bar\xi$ is $\cG$-complete it follows that
$\bar \Phi(\cdot, \tilde{x})= \widetilde{\Phi}(\cdot,
\tilde{x})$, for all $\tilde{x} \in 
(\{0\}\times \R^{m-l})^{\widetilde{}}_c$. Therefore $\bar \Phi
=\widetilde{\Phi}$ on $(\{0\}\times \R^{m-l})^{\widetilde{}}_c$
and $\bar \Phi(\eta, \tilde{x})\in S_{\bar F}$,
for all $\eta$ and all $\tilde{x} \in S_{\bar F}$,
i.e. $\bar \Phi$ is a symmetry of $\bar F=0$.

Alternatively, let us assume that (ii) obtains. 
We have to show that $\mbox{pr}_1 \circ \bar \Phi
(\tilde{\eta}, \tilde{x})=0$ for all $\tilde{\eta}
\in \widetilde{\R}_c$ and $\tilde{x} \in S_F$. Let $1\leq k\leq l$. Then
for representatives as above, $\bar\Phi_{k\eps}(0, x_\eps)=0$ and 
since $V$ is a box, $(\sigma
\bar\Phi_{\eps}'(\tau, x_\eps), \bar\Phi_{\eps}''(\tau, x_\eps)) \in V$
for $\sigma \in [0,1]$ and $\tau\in (-\eta_0,\eta_0)$. Therefore,  
\beast
\bar\Phi_{k\eps}(\eta, x_\eps) & = & \int_0^\eta
\frac{d}{d\tau}\bar\Phi_{k\eps}(\tau, x_\eps) \, d\tau 
=\int_0^\eta \bar\xi_{k\eps}(\bar\Phi_\eps(\tau, x_\eps)) \, d\tau \\
& = & \int_0^\eta 
\Big(\bar\xi_{k\eps}(\underbrace{\bar\Phi_{1\eps}(\tau, x_\eps),
\dots, \bar\Phi_{l\eps}(\tau, x_\eps)}_{=:\bar\Phi_{\eps}'(\tau, 
x_\eps)},
\underbrace{\bar\Phi_{{l+1}\eps}(\tau, x_\eps), \dots,
\bar\Phi_{m\eps}(\tau, x_\eps)}_{=:\bar\Phi_{\eps}''(\tau, x_\eps)}) \\
&   & \qquad - \bar\xi_{k\eps}(0, \bar\Phi_{\eps}''(\tau, 
x_\eps))\Big)\, d\tau \\
& = & \int_0^\eta \int_0^1 \frac{d}{d\sigma}\bar\xi_{k\eps}(\sigma
\bar\Phi_{\eps}'(\tau, x_\eps), \bar\Phi_{\eps}''(\tau, x_\eps))\, 
d\sigma \, d\tau \\
& = & \int_0^\eta \int_0^1 \sum_{j=1}^l D_j\bar\xi_{k\eps}(\sigma
\bar\Phi_{\eps}'(\tau, x_\eps), \bar\Phi_{\eps}''(\tau, x_\eps))\cdot 
\bar\Phi_{j\eps}(\tau, x_\eps)
\, d\sigma \, d\tau.
\eeast
Since $\bar\Phi$ is c-bounded and $D\bar\xi$ is locally of 
$L^\infty$-$\log$
type, the claim therefore follows by applying Gronwall's inequality.
\ep
\brem We list some sufficient conditions for the respective 
assumptions of the above theorem:
\begin{itemize}
\item[(i)] In case $M$ is a Riemannian manifold with Riemannian metric 
$h$
(e.g., a submanifold of $\R^n$ with the induced metric) it suffices to 
assume that $\xi$ and $P\psi$ are globally bounded with respect to $h$
for each differential operator $P$ of first order. 
\item[(ii)] To secure this condition it suffices to suppose that $P\xi$
is locally bounded for each differential operator $P$ of order $\le 1$ 
and
that $P\psi$ is locally bounded for each differential operator $P$
of order $\le 2$.
\end{itemize}
\ethi

\bexs \label{chartex} 
In certain algebraically special cases 
a global chart $\psi$ as in Theorem \ref{glavna} can immediately
be read off: 
\begin{itemize}
        \item[(i)]
 Suppose that (after a possible renumbering of the coordinates)
$F\in \gs(\R^n)^l$ is given in the form  
$$
\begin{array}{lcl}
F_1(x_1,\dots,x_n) &=& x_1 - f_1(x_2,\dots,x_n) \\  
F_2(x_1,\dots,x_n) &=& x_2 - f_2(x_3,\dots,x_n) \\  
\qquad\quad\dots && \qquad\quad\dots\\
F_l(x_1,\dots,x_n) &=& x_l - f_l(x_{l+1},\dots,x_n)
\end{array}
$$
with $f_i \in \gs[\R^{n-i},\R]$ for $1\le i \le l$. Then
$$
\psi(y_1,\dots,y_n) = (F_1(y_1,\dots,y_n),\dots,F_l(y_1,\dots,y_n),
y_{l+1},\dots,y_n)
$$
and writing $\psi^{-1}(x_1,\dots,x_n)=(y_1,\dots,y_n)$, $\psi^{-1}$
is determined recursively by $y_i=x_i$ ($l< i\le n$) and 
$y_i = x_i + f_i(y_{i+1},\dots,y_n)$ for $i \le l$. Since composition
of c-bounded generalized functions can be carried out unrestrictedly
(\cite{gfvm2}, Th.\ 3.6), $\psi$ is a global generalized chart. For 
$l=1$
we obtain \cite{symm}, Th.\ 4.7 (formulated there in the 
$\gs_\tau$-setting) 
as a special case.
\item[(ii)] 
$A \in \widetilde{\R}_c^{n\times n}$ is an invertible matrix of generalized 
numbers
if and only if $\det(A)$ is strictly nonzero in $\widetilde \R$ (cf.\ 
\cite{book},
Th.\ 1.2.38 and Lemma 1.2.41). If, in addition, $\psi: x \mapsto 
A^{-1}\cdot x$
is c-bounded then it is a global chart for the map
$F\in \gs(\R^n)^l$, $F(x)=\mbox{pr}_{\R^n\rightarrow\R^l}(A\cdot x)$.
As a concrete example one may take for $A$ a generalized rotation, 
i.e., an element 
of the special orthogonal group $\mathsf{SO}(n,\widetilde{\R})$ over the 
ring $\widetilde{\R}$ of
generalized numbers (cf.\ \cite{ObCont, MoRot} and section \ref{rotsect}). 
\end{itemize}
\ethi

 \section{Differential Equations} \label{infcrit}
Based on the previous section, it is possible to derive a theory of 
symmetry groups of differential equations in the space of 
c-bounded generalized functions. This development largely parallels
the one presented in \cite{symm}, Sec.\ 4.2, though with the additional
benefit of being formulated in a global setting. Therefore we only 
point
out the technical differences and omit proofs which are analogous
to the $\gs_\tau$-setting used there.
 \bd
 A generalized group action $\Phi \in \cG[\R \times \R^{p+q},
 \R^{p+q}]$ is called projectable if
 \beq
 \Phi(\eta, (x,u))=(\Xi_\eta(x), \Psi_\eta(x,u)),
 \eeq
 where $\Xi\in \cG[\R \times \rp, \rp]$ and $\Psi \in
 \cG[\R \times \R^{p+q}, \rqq]$.
 \ethi

The group properties
 \beq \nonumber
 \begin{split}
 \Xi_{\eta_1 + \eta_2}(x) & = \Xi_{\eta_1}(\Xi_{\eta_2}(x))\\
 \Psi_{\eta_1 + \eta_2}(x,u) & = \Psi_{\eta_1}(\Xi_{\eta_2}(x),
 \Psi_{\eta_2}(x,u))
 \end{split}
 \eeq
are to be understood as equations in $\gs[\R^2\times\R^p,\R^p]$ and
$\gs[\R^2\times\R^{p+q},\R^q]$, respectively.
\cite{gfvm2}, Th.\ 3.5 shows that any element $u$ of $\gs[M,N]$ is
uniquely determined by its graph $\Gamma_u$. We have
 \bp
 Let $u\in \cG[\rp, \rqq]$ and let $\Phi$ be a projectable
 generalized group action on $\rp \times \rqq$. Then
 $\Phi_\eta(\Gamma_u)=\Gamma_{\Phi_\eta (u)}$ in
 $\widetilde{\R}^{p+q}_c$ for each $\eta \in \rtilda_c$,
 where $\Phi_\eta(u)$ denotes the element
 $$
 x \mapsto \Psi_\eta (\Xi_{-\eta}(x), u\circ \Xi_{-\eta}(x))
 \in \cG[\rp, \rqq].
 $$
 \ethi
 \bp \label{prop system}
 Consider a system of PDEs
 \beq \label{system}
 \Delta_\nu (x, \un)=0 \qquad\qquad (1\leq \nu \leq l)
 \eeq
 in $\cG[\rp, \rqq]$, where $\Delta \in \cG[(\rp \times
 \rqq)^{(n)}, \rl]$. Set
 $$
 S_\Delta = \{\tilde{z} \in (\rtilda^p_c)^{(n)}:\
 \Delta_\nu(\tilde{z})=0,\ 1\leq \nu \leq l \}.
 $$
 Then $u\in \cG[\rp, \rqq]$ is a solution of the system
 if and only if $\Gamma_{{\rm pr}^{(n)} u}\subseteq S_\Delta$.
 \ethi

 Prolongations of generalized group actions are constructed
 as in the classical theory: Let
 $\Phi$ be a projectable generalized group action on
 $\rp \times \rqq$, $z\in (\rp \times \rqq)^{(n)}$ and
 choose a function $h\in \cbesk(\rp, \rqq)$ such that
 $(z_1, \dots , z_p, \prn h(z_1, \dots , z_p))=z$. The
 $n$-th prolongation of $\Phi$ is defined as
 $$
 \prn \Phi(\eta, z):= (\Xi_\eta(z_1, \dots , z_p),
 \prn(\Phi_\eta(h))(\Xi_\eta(z_1,\dots ,z_p))).
 $$
 By \cite{book}, 3.2.59 it follows that $\prn \Phi\in
 \cG[\R\times(\R^{p+q})^{(n)}, (\R^{p+q})^{(n)}]$.
 As in \cite{symm}, Lemma 4.12 and Prop.\ 4.13 it is seen that this 
definition 
does not depend on the particular choice of $h$ and that $\prn \Phi$ is 
a generalized group action on $(\rp \times \rqq)^{(n)}$.

\bp \label{algsymmdgsymm}
Let $\Phi$ be a projectable generalized group action
on $\rp \times \rqq$ such that $\prn \Phi$ is a
symmetry group of the algebraic equation $\Delta(z)=0$.
Then $\Phi$ is a symmetry group of (\ref{system}). 
\ethi

\bd
Let $\xi$ be a $\cG$-complete generalized vector field. The $n$-th
prolongation of $\xi$ is the infinitesimal generator of the $n$-th
prolongation of the generalized group action $\Phi$ corresponding
to $\xi$:
$$
\prn \xi|_z = \deta \prn \Phi_\eta(z).
$$
If $\prn \xi$ is $\gs$-complete, then both $\xi$ and $\Phi$ are 
called $\cG$-$n$-complete.
\ethi

\bt
Let \label{mainth}
\beq \label{sistem}
\Delta_\nu (x, \prn u)=0 \qquad\qquad (1\leq \nu \leq l)
\eeq
be a system of partial differential equations with 
$\Delta\in \cG(\rp)^l$. Let $\Phi$ be a generalized group
action on $\rp \times \rqq$ with infinitesimal generator $\xi$ and
suppose that $\Delta$ and $\prn \Phi$ satisfy the assumtions of Th.\ 
\ref{glavna}.
If
 $$ 
\prn \xi(\Delta)(\tilde{z})=0 \qquad\quad \forall \tilde{z}
\in (\rtilda^p_c \times \rtilda^q_c)^{(n)} \quad
\mbox{ with } \Delta(\tilde{z})=0,
 $$
then $\Phi$ is a symmetry group of (\ref{sistem}).  
\ethi
\pr Immediate from Th.\ \ref{glavna} and Prop.\ \ref{algsymmdgsymm}.\ep

As in \cite{symm}, Th.\ 4.17 we may now conclude that the classical 
algorithm for determining symmetries of a given system of differential
equations carries over to the generalized setting: make an ansatz for
the infinitesimal generators, calculate the prolongations according
to the classical formulas (cf.\ \cite{Olv}, Th.\ 2.36) and then apply
Th.\ \ref{mainth} to derive a system of determining equations in the
space of c-bounded Colombeau functions. 
Solutions of this system verifying the conditions of Th.\ \ref{mainth}
yield generalized symmetries of (\ref{sistem}). For recent applications 
to
weak solutions of systems of conservation laws we refer to \cite{SK}.

\section{Group invariant generalized functions} \label{rotsect}

In this final section we analyze the notion of invariance of Colombeau 
generalized functions
under generalized group actions. As in classical analysis and 
distribution theory this concept 
plays an important role with respect to applications (cf.\ the 
calculation of group invariant 
fundamental solutions in $\D'$ resp.\ $\gs$ in \cite{ga, wi, ObCont}).

We shall need the fact that composition of Colombeau
generalized functions and c-bounded generalized functions is 
always well-defined: 
\blem
Let $u=[(u_\eps)_\eps] \in \gs[M,N]$, $v=[(v_\eps)_\eps]\in \gs(N)$.
Then $v\circ u := [(v_\eps\circ u_\eps)_\eps]$ is a well-defined
element of $\gs(M)$.
\ethi
\pr To show that $(v_\eps\circ u_\eps)_\eps\in \esm(M)$, let $K\comp V$
for some chart $(V,\varphi)$ in $M$.
Since $v$ is c-bounded, there exist $K'\comp N$ and $\eps_0>0$ such 
that 
$v_\eps(K)\subseteq K'$ for all $\eps<\eps_0$. Without loss of 
generality
we may assume that $K'$ is contained in some chart $(W,\psi)$ of $N$.
Then the moderateness estimates for $v_\eps\circ u_\eps
= (v_\eps\circ \psi^{-1})\circ (\psi\circ u_\eps)$ on $K$ follow from 
the
chain rule and the respective estimates for $(u_\eps)_\eps$ and 
$(v_\eps)_\eps$. Suppose now that $[(u_\eps)_\eps] = [(u_\eps')_\eps]$
in $\gs[M,N]$ and let $\tilde{x} = [(x_\eps)_\eps] \in \widetilde{M}_c$.
Then by \cite{book}, Prop.\ 3.2.56, $[(u_\eps(x_\eps))] = 
[(u_\eps'(x_\eps))]$ in $\widetilde{M}_c$ and hence $[(v_\eps\circ 
u_\eps(x_\eps))] = [(v_\eps\circ u_\eps'(x_\eps))]$ in $\widetilde{M}_c$ 
by \cite{book}, Prop.\ 3.2.7. By the same result, $[(v_\eps\circ 
u_\eps)_\eps]
= [(v_\eps\circ u_\eps')_\eps]$ in $\gs(M)$. Finally, if $(v_\eps)_\eps
\in \ns(N)$ it is immediate from the c-boundedness of $(u_\eps)_\eps$
that $(v_\eps\circ u_\eps)_\eps \in \ns(M)$. Hence $v\circ u$ is 
well-defined,
as claimed. 
\ep
In particular, for $\Phi$ a generalized group action on $M$ and $f\in 
\gs(M)$,
it follows that $f\circ \Phi$ is a well-defined element of 
$\gs(\R\times M)$.
\bd
Let $\Phi$ be a generalized group action on $M$ and let $f\in 
\gs(M)$. $f$ is called invariant under $\Phi$ if $f\circ \Phi
= f\circ \pi_2$ in $\gs(\R\times M)$
(with $\pi_2: \R\times M \to M$, $\pi_2(\eta,x)=x$).
\ethi

 By the point value characterization of generalized functions
 (cf.\ \cite{book}, Th.\ 3.2.8) the above condition can equivalently
 be stated as follows:
 $$
 f(\Phi(\tilde{\eta},\tilde{x})) = f(\tilde{x}) \qquad \forall
 \tilde{\eta} \in \widetilde{\R}_c\
 \forall \tilde{x} \in \widetilde{M}_c\,.
 $$

 \bp
 Let $f\in\cG(M)$ and let $\Phi$ be a generalized group
 action on $M$ with infinitesimal generator $\xi$. Then
 the following statements are equivalent:
 \begin{itemize}
 \item[(i)] $f$ is $\Phi$-invariant
 \item[(ii)] $\xi(f)=0$ in $\cG(M)$.
 \end{itemize}
 \ethi

 \pr
 (i)$\Rightarrow$(ii) Since $f$ is $\Phi$-invariant we have
 $
 0=\frac{d}{d\eta}{\vert}_{_{0}}(f(\Phi(\eta,x)))=
 \xi(f)|_{x}$ in $\cG(M)$.
 
 (ii)$\Rightarrow$(i) Conversely, let $\xi(f)=0$ in $\cG(M)$. Then
 $$
 \frac{d}{d\eta}f(\Phi(\eta,\tilde{x})) = \xi(f)|_{\Phi(\eta,
 \tilde{x})} =0 \qquad \mbox{in } \gs(\R) \ \forall\tilde{x} \in 
 \widetilde{M}_c.
 $$
 Therefore, for each $\tilde{x}$ the map $\eta \mapsto f(\Phi(\eta, 
 \tilde{x}))$
 is constant in $\gs(\R)$, so $f\circ\Phi = f\circ \pi_2$, again by 
 \cite{book},
 Th.\ 3.2.8.\ep

Invariance properties of Colombeau generalized functions under 
generalized group 
actions have first been studied in \cite{ObCont, MoRot}. In particular, 
the following basic 
result was derived (\cite{MoRot}, Th.\ 2, formulated there in the 
$\gs_\tau$-setting):
\bt \label{tranlth}
Let $u\in \gs(\R^n)$. The following are equivalent:
\begin{itemize}
\item[(i)] $u(\tilde{x}_1+\eta,\tilde{x}_2,\dots,\tilde{x}_n) =
 u(\tilde{x})$ for all $\tilde{x}
\in \widetilde{\R}^n_c, \eta \in \widetilde{\R}_c$.
\item[(ii)] $\partial_{x_1} u = 0$ in $\gs(\R^n)$.
\item[(iii)] $u$ has a representative $(u_\eps)_\eps$ such that 
$\partial_{x_1}u_\eps\equiv 0$
for all $\eps$.
\end{itemize}
\ethi
It remained an open question there whether (i)--(iii) is equivalent to 
\begin{itemize}
\item[(i')] $u(\tilde{x}_1+\eta,\tilde{x}_2,\dots,\tilde{x}_n) =
 u(\tilde{x})$ for all $\tilde{x}
\in \widetilde{\R}^n_c, \eta \in \R$.
\end{itemize}
i.e., whether standard translations suffice to characterize 
translational invariance of 
Colombeau generalized functions. Meanwhile, Pilipovi\'c, Scarpalezos and
Valmorin have provided two alternative 
proofs (based on a Baire argument resp.\ on the construction of a 
parametrix) which show
that this question can be answered affirmatively (\cite{PSV}). In 
what follows we shall make
use of this result to resolve a further open question raised in 
\cite{ObCont} in the context of
generalized rotations. 

Recall from \ref{chartex} that we denote by $\mathsf{SO}(n,\widetilde{\R})$ 
the space of generalized 
rotations. Rotational invariance of Colombeau functions has been 
characterized in \cite{ObCont,MoRot} and has
been employed there to provide a new method of calculating rotationally 
invariant fundamental solutions, e.g.,
of the Laplace equation. The main characterization result is as follows 
(\cite{ObCont}, Th.\ 4.2):
\bt \label{rotth}
Let $u\in \gs(\R^n)$. The following are equivalent:
\begin{itemize}
        \item[(i)] $u\circ A = u$ in $\gs(\R^n)$ for all $A\in 
\mathsf{SO}(n,\widetilde{\R})$.
        \item[(ii)] $ \xi  u = 0$ in $\gs(\R^n)$ for all infinitesimal 
generators of $\mathsf{SO}(n,\R)$.
        \item[(iii)] $u$ possesses a representative $(u_\eps)_\eps$ 
such that each $u_\eps$ is rotationally invariant.   
\end{itemize}
\ethi
The following result affirmatively answers an open question from 
\cite{MoRot}:
\bt \label{rotprop} 
Items (i)--(iii) in Theorem \ref{rotth} are equivalent with
\begin{itemize}
        \item[(i')] $u\circ A = u$ in $\gs(\R^n)$ for all $A\in 
\mathsf{SO}(n,\R)$,
\end{itemize}
i.e., standard rotations suffice to characterize rotational invariance 
of Colombeau generalized functions.
\ethi
\pr Obviously (i) implies (i'). To prove the converse we first treat 
the case $n=2$. 
Let $\widetilde{A}\in \mathsf{SO}(2,\widetilde{\R})$. Then by \cite{MoRot}, Sec.\ 
3, Lemma 1 there exists some 
$\tilde{\eta} \in \widetilde{\R}_c$ such that 
$$
\widetilde{A} = \left[ \left(
\begin{array}{cr}
\cos(\eta_\eps) & -\sin(\eta_\eps) \\
\sin(\eta_\eps) & \cos(\eta_\eps)
\end{array}
\right)_{\!\!\eps}\right]
$$
Given $\tilde{x}$, $\tilde{y} \in \widetilde{\R}_c$ we have to show that 
$u(\widetilde{A} \cdot (\tilde{x},\tilde{y})^t) = u(\tilde{x},\tilde{y})$
in $\widetilde{\R}$. We may write $(\tilde{x},\tilde{y}) = [(r_\eps 
\cos(\theta_\eps),r_\eps\sin(\theta_\eps))]$ for suitable
$r_\eps \ge 0$, $\theta_\eps$. Now set $v_\eps := \theta \mapsto 
u_\eps(r_\eps \cos(\theta),r_\eps\sin(\theta))$. Then
$v=[(v_\eps)_\eps]\in \gs(\R)$ and by assumption $v(\tilde{\theta} 
+\eta) = v(\tilde{\theta})$ in $\widetilde{\R}$ for all $\tilde{\theta}
\in \widetilde{\R}_c$ and all $\eta\in \R$. But then the equivalence of (i) 
and (i') in Th.\ \ref{tranlth} shows that $v$ is in fact
a generalized constant. This immediately gives the result in the 
2D-case.

In the general case $n\ge 2$ we verify (ii) of Th.\ \ref{rotth}. Let 
$1\le i < j \le n$ and let $ \xi  = x_i\partial_{x_j} - x_j 
\partial_{x_i}$
be an infinitesimal generator of $\mathsf{SO}(n,\R)$. Fix compactly 
supported generalized numbers $\tilde{x}_1$,\dots,$\tilde{x}_{i-1}$,
$\tilde{x}_{i+1}$,\dots, $\tilde{x}_{j-1}$,
 $\tilde{x}_{j+1}$,\dots,$\tilde{x}_{n}$ and consider the maps
$$
w_\eps: (x_i,x_j) \mapsto u_\eps(\tilde{x}_1,\dots,
\tilde{x}_{i-1},x_i,\dots,\tilde{x}_{j-1},x_j,\dots,\tilde{x}_n)\,.
$$
Then $w=[(w_\eps)_\eps] \in \gs(\R^2)$ and from our assumption it 
follows that $w\circ A = w$ in $\gs(\R^2)$ for all
$A\in \mathsf{SO}(\R^2)$. By what we have already proved in the 2D-case 
and Th.\ \ref{rotth} it follows that 
$ \xi  w = 0$ in $\gs(\R^2)$. Hence from the point value 
characterization of Colombeau generalized functions it follows
that $ \xi  u=0$ in $\gs(\R^n)$ for each $ \xi  \in\mathsf{SO}(\R^n)$, 
as claimed. 
\ep
\brem
Let $\tilde{a}$ be a strictly nonzero (i.e., invertible, cf.\ 
\cite{book}, Th.\ 1.2.38) generalized number and consider
the generalized vector field $\xi = \tilde{a}(y\partial_x - x 
\partial_y)$ on $\R^2$. Then $\psi=[(\psi_\eps)_\eps]$
with
$$
\psi_\eps: (r,\theta) \mapsto (r \cos(a_\eps\theta),r \sin(a_\eps 
\theta))
$$
is a generalized chart in $\gs[\R^+\times (0,2\pi),\R^2\setminus 
(\R_0^+\times \{0\})]$. Moreover, the pullback
$\psi^*\xi$ of $\xi$ under $\psi$ is the smooth vector field 
$\frac{\partial}{\partial \theta}$. This provides a
simple case of ``straightening out'' a (strictly) nonzero generalized 
vector field. In the case of standard polar
coordinates ($\tilde{a} = 1$) $\psi$ allows to directly transform 
standard generators of $\mathsf{SO}(2,\R)$ to
translations, albeit only on $\R^2\setminus\{0\}$. However, there exist 
elements of $\gs(\R^2)$ with
support $\{0\}$ which are not rotationally invariant: choose some 
$\vphi\in \D(\R^2)$ whose support is not
rotationally invariant and set $u=[(\vphi(\frac{.}{\eps}))_\eps]$. 
Therefore, the above argument does not yield an
alternative proof of Prop.\ \ref{rotprop} (by reducing it to the 
translation setting of Th.\ \ref{tranlth}), 
since, contrary to the smooth setting, 
rotational invariance on $\R^2\setminus\{0\}$ is not equivalent to 
rotational invariance on $\R^2$ for 
Colombeau functions. (The situation for $\D'(\R^2)$ is similar: for 
example, $\partial_1\delta$ is a distribution
supported in $\{0\}$ which is not rotationally invariant).

Nevertheless, generalized charts induced by matrix transformations as 
above and the related question of 
straightening out infinitesimal generators of matrix groups over the 
ring of generalized numbers are
likely to play an important role in a further analysis of group 
invariant generalized functions. 
They should also provide valuable test cases for the development of 
inverse function theorems in the
Colombeau setting (\cite{EG}).
\ethi


\begin{thebibliography}{10}

\bibitem{ga}
{Berest, Yu. Yu.}
\newblock Group analysis of linear differential equations in distributions and
  the construction of fundamental solutions.
\newblock {\em Diff.\ Equ.}, {\bf 29}(11):1700--1711, 1993.

\bibitem{wi}
{Berest, Yu. Yu.}
\newblock Weak invariants of local groups of transformations.
\newblock {\em Diff.\ Equ.}, {\bf 29}(10):1561--1567, 1993.

\bibitem{bi}
{Berest, Yu. Yu., Ibragimov, N. H.}
\newblock Group theoretic determination of fundamental solutions.
\newblock {\em Lie Groups Appl.}, {\bf 1}(2):65--80, 1994.

\bibitem{BK}
{Bluman, G.~W., Kumei, S.}
\newblock {\em Symmetries and Differential Equations}.
\newblock Springer, New York, 1989.

\bibitem{c1}
{Colombeau, J.~F.}
\newblock {\em New Generalized Functions and Multiplication of Distributions}.
\newblock North Holland, Amsterdam, 1984.

\bibitem{c2}
{Colombeau, J.~F.}
\newblock {\em Elementary Introduction to New Generalized Functions}.
\newblock North Holland, Amsterdam, 1985.

\bibitem{DKP}
{Djapi{\'c}, N., Kunzinger, M., Pilipovi{\'c}, S.}
\newblock Symmetry group analysis of weak solutions.
\newblock {\em Proc. London Math. Soc. (3)}, 84(3):686--710, 2002.

\bibitem{EG}
{Erlacher, E., Grosser, M.}
\newblock An inverse function theorem for generalized functions.
\newblock {\em Preprint}, 2004.

\bibitem{found}
{Grosser, M., Farkas, E., Kunzinger, M., Steinbauer, R.}
\newblock On the foundations of nonlinear generalized functions {I}, {II}.
\newblock {\em Mem. Amer. Math. Soc.}, {\bf 153}(729), 2001.

\bibitem{book}
{Grosser, M., Kunzinger, M., Oberguggenberger, M., Steinbauer, R.}
\newblock {\em Geometric Theory of Generalized Functions}, volume 537 of {\em
  Mathematics and its Applications 537}.
\newblock Kluwer Academic Publishers, Dordrecht, 2001.

\bibitem{amp}
{Ibragimov, N. H.}
\newblock Group theoretical treatment of fundamental solutions.
\newblock In {\em Analysis, Manifolds and Physics}. Kluwer, Dordrecht, 1992.

\bibitem{SK}
{Konjik, S.}
\newblock Symmetries of conservation laws.
\newblock {\em Publ. Inst. Math. (Beograd) (N.S.)}, to appear.

\bibitem{gfvm}
{Kunzinger, M.}
\newblock Generalized functions valued in a smooth manifold.
\newblock {\em Monatsh.\ Math.}, {\bf 137}:31--49, 2002.

\bibitem{symm}
{Kunzinger, M., Oberguggenberger, M.}
\newblock Group analysis of differential equations and generalized functions.
\newblock {\em SIAM J. Math. Anal.}, {\bf 31}(6):1192--1213, 2000.

\bibitem{flows}
{Kunzinger, M., Oberguggenberger, M., Steinbauer, R., Vickers, J.}
\newblock Generalized flows and singular {ODE}s on differentiable manifolds.
\newblock {\em Acta Appl. Math.}, {\bf 80}:221--241, 2004.

\bibitem{ndg}
{Kunzinger, M., Steinbauer, R.}
\newblock Foundations of a nonlinear distributional geometry.
\newblock {\em Acta Appl. Math.}, 71:179--206, 2002.

\bibitem{gprg}
{Kunzinger, M., Steinbauer, R.}
\newblock Generalized {pseudo-}{R}iemannian geometry.
\newblock {\em Trans. Amer. Math. Soc.}, 354(10):4179--4199, 2002.

\bibitem{gfvm2}
{Kunzinger, M., Steinbauer, R., Vickers, J.}
\newblock Intrinsic characterization of manifold-valued generalized functions.
\newblock {\em Proc.\ London Math.\ Soc.}, {\bf 87}(2):451--470, 2003.

\bibitem{Me}
{Meth\'ee, P. D.}
\newblock Sur les distributions invariantes dans le groupe des rotations de
  {L}orentz.
\newblock {\em Comment.~Math.~Helv.}, {\bf 28}:224--269, 1954.

\bibitem{MOBook}
{Oberguggenberger, M.}
\newblock {\em Multiplication of Distributions and Applications to Partial
  Differential Equations}, {\em Pitman Research Notes in
  Mathematics} {\bf 259}.
\newblock Longman, Harlow, U.K., 1992.

\bibitem{ObCont}
{Oberguggenberger, M.}
\newblock Symmetry groups, nonlinear partial differential equations, and
  generalized functions.
\newblock In {Lesley, J.\ A., Robart, T.}, (Ed.), {\em Geometrical Study of
  Differential Equations}, {\em Contemporary Mathematics} {\bf 285},
  101--110. Amer. Math. Soc., 2001.

\bibitem{MoRot}
{Oberguggenberger, M.}
\newblock Rotationally invariant {C}olombeau functions.
\newblock In {Delcroix, A., Hasler, M., Marti, J.-A., Valmorin, V.}, editor,
  {\em {N}onlinear {A}lgebraic {A}nalysis and {A}pplications}, 227--236.
  Cambridge Scientific Publishers, 2004.

\bibitem{point}
{Oberguggenberger, M., Kunzinger, M.}
\newblock Characterization of {C}olombeau generalized functions by their
  pointvalues.
\newblock {\em Math. Nachr.}, {\bf 203}:147--157, 1999.

\bibitem{OR}
{Oberguggenberger, M., Rosinger, E.~E.}
\newblock {\em Solution of Continuous Nonlinear PDEs through Order Completion}.
\newblock North Holland, Amsterdam, 1994.

\bibitem{Olv}
{Olver, P.~J.}
\newblock {\em Applications of Lie Groups to Differential Equations},
  volume~{\bf 107} of {\em Graduate Texts in Mathematics}.
\newblock Springer, New York, second edition, 1993.

\bibitem{PSV}
{Pilipovi\'c, S., Scarpalezos, D., Valmorin, V.}
\newblock Equalities in algebras of generalized functions.
\newblock {\em Preprint}, 2004.

\bibitem{RW}
{Rosinger, E.~E., Walus, Y.~E.}
\newblock Group invariance of generalized solutions obtained through the
  algebraic method.
\newblock {\em Nonlinearity}, {\bf 7}:837--859, 1994.

\bibitem{RW1}
{Rosinger, E.~E., Walus, Y.~E.}
\newblock Group invariance of global generalized solutions of nonlinear
  {P}{D}{E}s in nowhere dense algebras.
\newblock {\em Lie Groups Appl.}, {\bf 1}:216--225, 1994.

\bibitem{S}
{Szmydt, Z.}
\newblock On homogeneous rotation invariant distributions and the laplace
  operator.
\newblock {\em Ann.~Pol.~Math.}, {\bf 6}:249--259, 1979.

\bibitem{SZ}
{Szmydt, Z., Ziemian, B.}
\newblock Invariant fundamental solutions of the wave operator.
\newblock {\em Demonstr.~Math.}, {\bf 19}:371--386, 1986.

\bibitem{Zie}
{Ziemian, B.}
\newblock On distributions invariant with respect to some linear
  transformations.
\newblock {\em Ann. Pol. Math.}, {\bf 36}:261--276, 1979.

\end{thebibliography}

\end{document}